\documentclass[11pt]{article}

\usepackage[matrix,arrow,curve,cmtip]{xy}
\usepackage{amssymb}
\usepackage{latexsym}
\usepackage{theorem}
\usepackage{titlesec}
\usepackage{graphicx}

\usepackage{times}

\usepackage[a4paper,text={128mm,195mm}]{geometry}

\titleformat{\section}[hang]%
{\bfseries\large}{\thesection.}{1ex}{}%

\titleformat{\subsection}[hang]%
{\bfseries}{\thesubsection}{1ex}{}%

\def\to{\mbox{$\xymatrix@1@C=5mm{\ar@{->}[r]&}$}}
\def\tto{\mbox{$\xymatrix@1@C=5mm{\ar@{=>}[r]&}$}}
\def\halfcirc{\begin{picture}(0,0)\put(0,2){\oval(4,4)[l]}\end{picture}}
\def\incl{\mbox{$\xymatrix@1@C=5mm{\ar@{->}[r]|<{\halfcirc}&}$}}
\def\distsign{\begin{picture}(0,0)\put(0,0){\circle{4}}\end{picture}}
\def\dist{\mbox{$\xymatrix@1@C=5mm{\ar@{->}[r]|{\distsign}&}$}}
\def\biar{\mbox{$\xymatrix@1@C=5mm{\ar@<1.5mm>[r]\ar@<-0.5mm>[r]&}$}}
\def\bidist{\mbox{$\xymatrix@1@C=5mm{\ar@<1.5mm>[r]|{\distsign}\ar@<-0.5mm>[r]|{\distsign}&}$}}
\def\iso{\mbox{$\xymatrix@1@C=6mm{\ar@{->}[r]^{\sim}&}$}}
\def\inlineadj#1#2{\mbox{$\xymatrix@C=15mm{\ar@{}[r]|{\bot}\ar@<1mm>@/^2mm/[r]^{{#1}} & \ar@<1mm>@/^2mm/[l]^{{#2}}}$}}

\newtheorem{theorem}{Theorem}[section]
\newtheorem{lemma}[theorem]{Lemma}
\newtheorem{definition}[theorem]{Definition} 
\newtheorem{proposition}[theorem]{Proposition}
\newtheorem{corollary}[theorem]{Corollary}
{\theorembodyfont{\upshape}\newtheorem{example}[theorem]{Example}}
\newcommand{\proof}{\noindent {\it Proof\ }: }
\def\endofproof{$\mbox{ }\hfill\Box$\par\vspace{1.8mm}\noindent}

\def\eqref#1{(\ref{#1})}

\def\:{\colon}
\def\impl{\Longrightarrow}
\def\2{{\bf 2}}

\def\op{^{\sf op}}
\def\Sup{{\sf Sup}}
\def\Dist{{\sf Dist}}
\def\Cat{{\sf Cat}}

\def\Q{{\cal Q}}

\def\C{{\cal C}}
\def\P{{\cal P}}
\def\V{{\cal V}}
\def\H{{\cal H}}
\def\K{{\cal K}}
\def\T{{\cal T}}
\def\colim{\mathop{\rm colim}}
\def\lim{\mathop{\rm lim}}
\def\bbA{\mathbb{A}}
\def\bbB{\mathbb{B}}
\def\bbC{\mathbb{C}}

\def\bbI{\mathbb{I}}

\def\tensor{\otimes}

\def\<{\langle}
\def\>{\rangle}

\begin{document}





\title{\Large `Hausdorff distance' via conical cocompletion}
\author{Isar Stubbe}
\date{\small March 16, 2009, revised November 5, 2009}

\maketitle






\begin{abstract}
In the context of quantaloid-enriched categories, we explain how each saturated class of weights defines, and is defined by, an essentially unique full sub-KZ-doctrine of the free cocompletion KZ-doctrine. The KZ-doctrines which arise as full sub-KZ-doctrines of the free cocompletion, are characterised by two simple ``fully faithfulness'' conditions. Conical weights form a saturated class, and the corresponding KZ-doctrine is precisely (the generalisation to quantaloid-enriched categories of) the Hausdorff doctrine of [Akhvlediani {\it et al.}, 2009].
\end{abstract}


\section{Introduction}

At the meeting on ``Categories in Algebra, Geometry and Logic'' honouring Francis Borceux and Dominique Bourn in Brussels on 10--11 October 2008, Walter Tholen gave a talk entitled ``On the categorical meaning of Hausdorff and Gromov distances'', reporting on joint work with Andrei Akhvlediani and Maria Manuel Clementino [2009]. The term `Hausdorff distance' in his title refers to the following construction: if $(X,d)$ is a metric space and $S,T\subseteq X$, then
$$\delta(S,T):=\bigvee_{s\in S}\bigwedge_{t\in T}d(s,t)$$
defines a (generalised) metric on the set of subsets of $X$. But Bill Lawvere [1973] showed that metric spaces are examples of enriched categories, so one can aim at suitably generalising this `Hausdorff distance'. Tholen and his co-workers achieved this for {\em categories enriched in a commutative quantale $\V$}. In particular they devise a KZ-doctrine on the category of $\V$-categories, whose algebras -- in the case of metric spaces -- are exactly the sets of subsets of metric spaces, equipped with the Hausdorff distance.

We shall argue that the notion of Hausdorff distance can be developed for {\em quant\-aloid-enriched categories} too, using {\em enriched colimits} as main tool. In fact, very much in line with the work of [Albert and Kelly, 1988; Kelly and Schmitt, 2005; Schmitt, 2006] on cocompletions of categories enriched in a symmetric monoidal category and the work of [Kock, 1995] on the abstraction of cocompletion processes, we shall see that, for quantaloid-enriched categories, each saturated class of weights defines, and is defined by, an essentially unique KZ-doctrine. The KZ-doctrines that arise in this manner are the full sub-KZ-doctrines of the free cocompletion KZ-doctrine, and they can be characterised with two simple ``fully faithfulness'' conditions. As an application, we find that the conical weights form a saturated class and the corresponding KZ-doctrine is precisely (the generalisation to quantaloid-enriched categories of) the Hausdorff doctrine of [Akhvlediani {\it et al.}, 2009].

In this paper we do not speak of `Gromov distances', that other metric notion that Akhvlediani, Clementino and Tholen [2009] refer to. As they analyse, Gromov distance is necessarily built up from {\em symmetrised} Hausdorff distance; and because their base quantale $\V$ is commutative, they can indeed extend this notion too to $\V$-enriched categories. More generally however, symmetrisation for quantaloid-enriched categories makes sense when that quantaloid is involutive. Preliminary computations indicate that `Gromov distance' ought to exist on this level of generality, but quickly got too long to include them in this paper: so we intend to work this out in a sequel.

\section{Preliminaries}

\subsection{Quantaloids}

A {\bf quantaloid} is a category enriched in the monoidal category $\Sup$ of complete lattices (also called sup-lattices) and supremum preserving functions (sup-morphisms). A quantaloid with one object, i.e.\ a monoid in $\Sup$, is a {\bf quantale}. Standard references include [Rosenthal, 1996; Paseka and Rosicky, 2000] .

Viewing $\Q$ as a locally ordered category, the 2-categorical notion of {\em adjunction in $\Q$} refers to a pair of arrows, say $f\:A\to B$ and $g\:B\to A$, such that $1_A\leq g\circ f$ and $f\circ g\leq 1_B$ (in which case $f$ is left adjoint to $g$, and $g$ is right adjoint to $f$, denoted $f\dashv g$).

Given arrows
$$
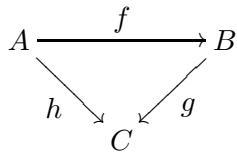
\xymatrix{A\ar[rr]^f\ar[dr]_h & & B\ar[dl]^g \\ & C}$$
in a quantaloid $\Q$, there are adjunctions between sup-lattices as follows:
$$\Q(B,C)\inlineadj{-\circ f}{\{f,-\}}Q(A,C),\quad \Q(A,B)\inlineadj{g\circ-}{[g,-]}\Q(A,C),$$ 
$$\Q(A,B)\inlineadj{\{-,h\}}{[-,h]}\Q(B,C)\op.$$
The arrow $[g,h]$ is called the {\bf lifting} of $h$ trough $g$, whereas $\{f,h\}$ is the {\bf extension} of $h$ through $f$. Of course, every left adjoint preserves suprema, and every right adjoint preserves infima. For later reference, we record some straightforward facts:
\begin{lemma}\label{a1}
If $g\:B\to C$ in a quantaloid has a right adjoint $g^*$, then $[g,h]=g^*\circ h$ and therefore $[g,-]$ also preserves suprema. Similarly, if $f\:A\to B$ has a left adjoint $f_!$ then $\{f,h\}=h\circ f_!$ and thus $\{f,-\}$ preserves suprema.
\end{lemma} 
\begin{lemma}\label{a2}
For any commutative diagram
$$\xymatrix{
A\ar[dr]_f & & B\ar[dl]_g\ar[dr]^h & & C\ar[dl]^i \\
 & D\ar[rr]_j & & E }$$
in a quantaloid, we have that $[i,h]\circ[g,f]\leq[i,j\circ f]$. If all these arrows are left adjoints, and $g$ moreover satisfies $g\circ g^*=1_D$, then $[i,h]\circ[g,f]=[i,j\circ f]$.
\end{lemma}
\begin{lemma}\label{a3.0}
If $f\:A\to B$ in a quantaloid has a right adjoint $f^*$ such that moreover $f^*\circ f=1_A$, then $[f\circ x, f\circ y]=[x,y]$ for any $x,y\:X\biar A$.
\end{lemma}

\subsection{Quantaloid-enriched categories}

From now on $\Q$ denotes a {\em small} quantaloid. Viewing $\Q$ as a (locally ordered) bicategory, it makes perfect sense to consider categories enriched in $\Q$. Bicategory-enriched categories were invented at the same time as bicategories by Jean B\'enabou [1967], and further developed by Ross Street [1981, 1983]. Bob Walters [1981] particularly used quantaloid-enriched categories in connection with sheaf theory. Here we shall stick to the notational conventions of [Stubbe, 2005], and refer to that paper for additional details, examples and references.

A {\bf $\Q$-category} $\bbA$ consists of a set of objects $\bbA_0$, a type function $t\:\bbA_0\to\Q_0$, and $\Q$-arrows $\bbA(a',a)\:ta\to ta'$; these must satisfy identity and composition axioms, namely:
$$1_{ta}\leq\bbA(a,a)\mbox{\ \ and \ \ }\bbA(a'',a')\circ\bbA(a',a)\leq\bbA(a'',a).$$
A {\bf $\Q$-functor} $F\:\bbA\to\bbB$ is a type-preserving object map $a\mapsto Fa$ satisfying the functoriality axiom:
$$\bbA(a',a)\leq\bbB(Fa',Fa).$$
And a {\bf $\Q$-distributor} $\Phi\:\bbA\dist\bbB$ is a matrix of $\Q$-arrows $\Phi(b,a)\:ta\to tb$, indexed by all couples of objects of $\bbA$ and $\bbB$, satisfying two action axioms:
$$\Phi(b,a')\circ\bbA(a',a)\leq\Phi(b,a)\mbox{\ \ and \ \ }\bbB(b,b')\circ\Phi(b',b)\leq\Phi(b,a).$$

Composition of functors is obvious; that of distributors is done with a ``matrix'' multiplication: the composite $\Psi\tensor\Phi\:\bbA\dist\bbC$ of $\Phi\:\bbA\dist\bbB$ and $\Psi\:\bbB\dist\bbC$ has as elements
$$(\Psi\tensor\Phi)(c,a)=\bigvee_{b\in\bbB_0}\Psi(c,b)\circ\Phi(b,a).$$
Moreover, the elementwise supremum of parallel distributors $(\Phi_i\:\bbA\dist\bbB)_{i\in I}$ gives a distributor $\bigvee_i\Phi_i\:\bbA\dist\bbB$, and it is easily checked that we obtain a (large) quantaloid $\Dist(\Q)$ of $\Q$-categories and distributors. Now $\Dist(\Q)$ is a 2-category, so we can speak of adjoint distributors. In fact, any functor $F\:\bbA\to\bbB$ determines an adjoint pair of distributors: 
\begin{equation}\label{a3.1.0}
\bbA\xymatrix@=14mm{\ar@{}[r]|{\perp}\ar@<1mm>@/^2mm/[r]^{\bbB(-,F-)}|{\distsign} & \ar@<1mm>@/^2mm/[l]^{\bbB(F-,-)}|{\distsign}}\bbB.
\end{equation}
Therefore we can sensibly order parallel functors $F,G\:\bbA\biar\bbB$ by putting $F\leq G$ whenever $\bbB(-,F-)\leq\bbB(-,G-)$ (or equivalently, $\bbB(G-,-)\leq\bbB(F-,-)$) in $\Dist(\Q)$. Doing so, we get a locally ordered category $\Cat(\Q)$ of $\Q$-categories and functors, together with a 2-functor
\begin{equation}\label{a3.1}
i\:\Cat(\Q)\to\Dist(\Q)\:\Big(F\:\bbA\to\bbB\Big)\mapsto\Big(\bbB(-,F-)\:\bbA\dist\bbB\Big).
\end{equation}
(The local order in $\Cat(\Q)$ need not be anti-symmetric, i.e.\ it is not a partial order but rather a preorder, which we prefer to call simply an order.)

This is the starting point for the theory of quantaloid-enriched categories, including such notions as:
\begin{itemize}
\item {\bf fully faithful functor}: an $F\:\bbA\to\bbB$ for which $\bbA(a',a)=\bbB(Fa',Fa)$, or alternatively, for which the unit of the adjunction in \eqref{a3.1.0} is an equality,
\item {\bf adjoint pair}: a pair $F\:\bbA\to\bbB$, $G\:\bbB\to\bbA$ for which $1_{\bbA}\leq G\circ F$ and also $F\circ G\leq 1_{\bbB}$, or alternatively, for which $\bbB(F-,-)=\bbA(-,G-)$,
\item {\bf equivalence}: an $F\:\bbA\to\bbB$ which are fully faithful and essentially surjective on objects, or alternatively, for which there exists a $G\:\bbB\to\bbA$ such that $1_{\bbA}\cong G\circ F$ and $F\circ G\cong 1_{\bbB}$,
\item {\bf left Kan extension}: given $F\:\bbA\to\bbB$ and $G\:\bbA\to\bbC$, the left Kan extension of $F$ through $G$, written $\<F,G\>\:\bbC\to\bbB$, is the smallest such functor satisfying $F\leq\<F,G\>\circ G$,
\end{itemize}
and so on. In the next subsection we shall recall the more elaborate notions of presheaves, weighted colimits and cocompletions.

\subsection{Presheaves and free cocompletion}

If $X$ is an object of $\Q$, then we write $*_X$ for the one-object $\Q$-category, whose single object $*$ is of type $X$, and whose single hom-arrow is $1_X$.

Given a $\Q$-category $\bbA$, we now define a new $\Q$-category $\P(\bbA)$ as follows:
\begin{itemize}
\item objects: $(\P(\bbA))_0=\{\phi\:*_X\dist\bbA\mid X\in\Q_0\}$,
\item types: $t(\phi)=X$ for $\phi\:*_X\dist\bbA$,
\item hom-arrows: $\P(\bbA)(\psi,\phi)=$ (single element of) the lifting $[\psi,\phi]$ in $\Dist(\Q)$.
\end{itemize}
Its objects are {\bf (contravariant) presheaves} on $\bbA$, and $\P(\bbA)$ itself is the {\bf presheaf category} on $\bbA$. 

The presheaf category $\P(\bbA)$ {\bf classifies distributors} with codomain $\bbA$: for any $\bbB$ there is a bijection between $\Dist(\Q)(\bbB,\bbA)$ and $\Cat(\Q)(\bbB,\P(\bbA))$, which associates to any distributor $\Phi\:\bbB\dist\bbA$ the functor $Y_{\Phi}\:\bbB\to\P(\bbA)\:b\mapsto\Phi(-,b)$, and conversely associates to any functor $F\:\bbB\to\P(\bbA)$ the distributor $\Phi_F\:\bbB\dist\bbA$ with elements $\Phi_F(a,b)=(Fb)(a)$. In particular is there a functor, $Y_{\bbA}\:\bbA\to\P(\bbA)$, that corresponds with the identity distributor $\bbA\:\bbA\dist\bbA$: the elements in the image of $Y_{\bbA}$ are the {\bf representable presheaves on $\bbA$}, that is to say, for each $a\in\bbA$ we have $\bbA(-,a)\:*_{ta}\dist\bbA$. Because such a representable presheaf is a left adjoint in $\Dist(\Q)$, with right adjoint $\bbA(a,-)$, we can verify that
$$\P(\bbA)(Y_{\bbA}(a),\phi)=[\bbA(-,a),\phi]=\bbA(a,-)\tensor\phi=\phi(a).$$
This result is known as {\bf Yoneda's Lemma}, and implies that $Y_{\bbA}\:\bbA\to\P(\bbA)$ is a fully faithful functor, called the {\bf Yoneda embedding} of $\bbA$ into $\P(\bbA)$.

By construction there is a 2-functor 
$$\P_0\:\Dist(\Q)\to\Cat(\Q)\:(\Phi\:\bbA\dist\bbB)\mapsto(\Phi\tensor-\:\P(\bbA)\to\P(\bbB)),$$
which is easily seen to preserve local suprema. Composing this with the one in \eqref{a3.1} we define two more 2-functors:
\begin{equation}\label{a4}
\begin{array}{c}
\xymatrix@=15mm{
\Dist(\Q)\ar[dr]^{\P_0}\ar@{.>}[r]^{\P_1} & \Dist(\Q) \\
\Cat(\Q)\ar[u]^i\ar@{.>}[r]_{\P} & \Cat(\Q)\ar[u]_i}
\end{array}
\end{equation}
In fact, $\P_1$ is a $\Sup$-functor (a.k.a.\ a homomorphism of quantaloids). Later on we shall encounter these functors again.

For a distributor $\Phi\:\bbA\dist\bbB$ and a functor $F\:\bbB\to\bbC$ between $\Q$-categories, the {\bf $\Phi$-weighted colimit of $F$} is a functor $K\:\bbA\to\bbC$ such that 
$[\Phi,\bbB(F-,-)]=\bbC(K-,-)$. Whenever a colimit exists, it is essentially unique; therefore the notation $\colim(\Phi,F)\:\bbA\to\bbC$ makes sense. These diagrams picture the situation:
$$\xymatrix@=15mm{
\bbB\ar[r]^F & \bbC \\
\bbA\ar[u]|{\distsign}^{\Phi}\ar@{.>}[ru]_{\colim(\Phi,F)}}
\hspace{1cm}
\xymatrix@=15mm{
\bbB & \bbC\ar[l]|{\distsign}_{\bbC(F-,-)}\ar@{.>}[dl]|{\distsign}^{[\Phi,\bbC(F-,-)]=\bbC(\colim(\Phi,F)-,-)} \\
\bbA\ar[u]|{\distsign}^{\Phi} & }$$
A functor $G\:\bbC\to\bbC'$ is said to {\bf preserve} $\colim(\Phi,F)$ if $G\circ\colim(\Phi,F)$ is the $\Phi$-weighted colimit of $G\circ F$. A $\Q$-category admitting all possible colimits, is {\bf cocomplete}, and a functor which preserves all colimits which exist in its domain, is {\bf cocontinuous}. (There are, of course, the dual notions of limit, completeness and continuity. We shall only use colimits in this paper, but it is a matter of fact that a $\Q$-category is complete if and only if it is cocomplete [Stubbe, 2005, Proposition 5.10].)

For two functors $F\:\bbA\to\bbB$ and $G\:\bbA\to\bbC$, we can consider the $\bbC(G-,-)$-weighted colimit of $F$. Whenever it exists, it is $\<F,G\>\:\bbC\to\bbB$, the left Kan extension of $F$ through $G$; but not every left Kan extension need to be such a colimit. Therefore we speak of a {\bf pointwise left Kan extension} in this case.

Any presheaf category $\P(\bbC)$ is cocomplete, as follows from its classifying property: given a distributor $\Phi\:\bbA\dist\bbB$ and a functor $F\:\bbB\to\P(\bbC)$, consider the unique distributor $\Phi_F\:\bbB\to\bbC$ corresponding with $F$; now in turn the composition $\Phi_F\tensor\Phi\:\bbA\dist\bbC$ corresponds with a unique functor $Y_{\Phi_F\tensor\Phi}\:\bbA\to\P(\bbC)$; the latter is $\colim(\Phi,F)$.

In fact, the 2-functor 
$$\P\:\Cat(\Q)\to\Cat(\Q)$$ 
is the {\bf Kock-Z\"oberlein-doctrine\footnote{A Kock-Z\"oberlein-doctrine (or KZ-doctrine, for short) $\T$ on a locally ordered category $\K$ is a 2-functor $\T\:\K\to\K$ for which there are a multiplication $\mu\:\T\circ\T\tto\T$ and a unit $\eta\:1_{\K}\tto\T$ making $(\T,\mu,\eta)$ a 2-monad, and satisfying moreover the ``KZ-inequation'': $\T(\eta_K)\leq\eta_{\T(K)}$ for all objects $K$ of $\K$. The notion was invented independently by Volker Z\"oberlein [1976] and Anders Kock [1972] in the more general setting of 2-categories. We refer to [Kock, 1995] for all details.} for free cocompletion}; the components of its multiplication $M\:\P\circ\P\tto\P$ and its unit $Y\:1_{\Cat(\Q)}\tto\P$ are
$$\colim(-,1_{\P(\bbC)})\:\P(\P(\bbC))\to\P(\bbC)\mbox{\ \ and \ \ }Y_{\bbC}\:\bbC\to\P(\bbC).$$
This means in particular that $(\P, M, Y)$ is a monad on $\Cat(\Q)$, and a $\Q$-category $\bbC$ is cocomplete if and only if it is a $\P$-algebra, if and only if $Y_{\bbC}\:\bbC\to\P(\bbC)$ admits a left adjoint in $\Cat(\Q)$.

\subsection{Full sub-KZ-doctrines of the free cocompletion doctrine}

The following observation will be useful in a later subsection.
\begin{proposition}\label{a9.1}
Suppose that $\T\:\Cat(\Q)\to\Cat(\Q)$ is a 2-functor and that 
$$\xy
\xymatrix@=15mm{
\Cat(\Q)\ar@<0mm>@/^5mm/[r]^{\P}\ar@<-0mm>@/_5mm/[r]_{\T} & \Cat(\Q) }
\POS(15,0)\drop{\rotatebox{90}{$\Longrightarrow$}}\POS(17,0)\drop{\varepsilon}
\endxy$$
is a 2-natural transformation, with all components $\varepsilon_{\bbA}\:\T(\bbA)\to\P(\bbA)$ fully faithful functors, such that there are (necessarily essentially unique) factorisations
$$\xymatrix@R=5mm@C=15mm{
\P\circ\P\ar@{=>}[r]^{ M} & \P  \\
 & & 1_{\Cat(\Q)}\ar@{=>}[ul]_{ Y}\ar@{:>}[dl]^{\eta} \\
\T\circ\T\ar@{=>}[uu]^{\varepsilon*\varepsilon}\ar@{:>}[r]_{\mu} & \T\ar@{=>}[uu]_{\varepsilon}}$$
Then $(\T,\mu,\eta)$ is a sub-2-monad of $(\P, M, Y)$, and is a KZ-doctrine. We call the pair $(\T,\varepsilon)$ a {\bf full sub-KZ-doctrine} of $\P$.
\end{proposition}
\proof
First note that, because each $\varepsilon_{\bbA}\:\T(\bbA)\to\P(\bbA)$ is fully faithful, for each $F,G\:\bbC\to\T(\bbA)$,
$$\varepsilon_{\bbA}\circ F\leq\varepsilon_{\bbA}\circ G\impl F\leq G,$$
thus in particular $\varepsilon_{\bbA}$ is (essentially) a monomorphism in $\Cat(\Q)$: if $\varepsilon_{\bbA}\circ F\cong\varepsilon_{\bbA}\circ G$ then $F\cong G$. Therefore we can regard $\varepsilon\:\T\tto\P$ as a subobject of the monoid $(\P, M, Y)$ in the monoidal category of endo-2-functors on $\Cat(\Q)$. The factorisations of $ M$ and $ Y$ then say precisely that $(\T,\mu,\eta)$ is a submonoid, i.e.\ a 2-monad on $\Cat(\Q)$ too.

But $\P\:\Cat(\Q)\to\Cat(\Q)$ maps fully faithful functors to fully faithful functors, as can be seen by applying Lemma \ref{a3.0} to the left adjoint $\bbB(-,F-)\:\bbA\dist\bbB$ in $\Dist(\Q)$, for any given fully faithful $F\:\bbA\to\bbB$. Therefore each 
$$(\varepsilon*\varepsilon)_{\bbA}\:\T(\T(\bbA))\to\P(\P(\bbA))$$ 
is fully faithful: for $(\varepsilon*\varepsilon)_{\bbA}=\P(\varepsilon_{\bbA})\circ\varepsilon_{\T\bbA}$ and by hypothesis both $\varepsilon_{\bbA}$ and $\varepsilon_{\T\bbA}$ are fully faithful. The commutative diagrams
$$\xymatrix@=15mm{
\P(\bbA)\ar[r]^{\P( Y_{\bbA})} & \P(\P(\bbA)) \\
 & \T(\P(\bbA))\ar[u]^{Y_{\T(\bbA)}} \\
\T(\bbA)\ar[uu]^{\varepsilon_{\bbA}}\ar[r]_{\T(\eta_{\bbA})}\ar[ru]^{\T(Y_{\bbA})} & \T(\T(\bbA))\ar[u]^(0.3){\T(\varepsilon_{\bbA})}\ar@<-1mm>@/_8mm/[uu]_{(\varepsilon*\varepsilon)_{\bbA}}}
\hspace{5mm}
\xymatrix@=15mm{
\P(\bbA)\ar[r]^{Y_{\P(\bbA)}}\ar[dr]_{\eta_{\P\bbA}} & \P(\P(\bbA)) \\
 & \T(\P(\bbA))\ar[u]^{Y_{\T(\bbA)}} \\
\T(\bbA)\ar[uu]^{\varepsilon_{\bbA}}\ar[r]_{\eta_{\T(\bbA)}} & \T(\T(\bbA))\ar[u]^{\T(\varepsilon_{\bbA})}\ar@<-1mm>@/_8mm/[uu]_{(\varepsilon*\varepsilon)_{\bbA}}}$$
thus imply, together with the KZ-inequation for $\P$, the KZ-inequation for $\T$.
\endofproof

Some remarks can be made about the previous Proposition. Firstly, about the fully faithfulness of the components of $\varepsilon\:\T\tto\P$. In any locally ordered category $\K$ one defines an arrow $f\:A\to B$ to be {\em representably fully faithful} when, for any object $X$ of $\K$, the order-preserving function
$$\K(f,-)\:\K(X,A)\to\K(X,B)\:x\mapsto f\circ x$$
is order-reflecting -- that is to say, $\K(f,-)$ is a fully faithful functor between ordered sets viewed as categories -- and therefore $f$ is also essentially a monomorphism in $\K$. But the converse need not hold, and indeed does not hold in $\K=\Cat(\Q)$: not every monomorphism in $\Cat(\Q)$ is representably fully faithful, and not every representably fully faithful functor is fully faithful. Because the 2-functor $\P\:\Cat(\Q)\to\Cat(\Q)$ preserves representable fully faithfulness as well, the above Proposition still holds (with the same proof) when the components of $\varepsilon\:\T\tto\P$ are merely {\em representably} fully faithful; and in that case it might be natural to say that $\T$ is a ``sub-KZ-doctrine'' of $\P$. But for our purposes later on, the interesting notion is that of {\em full} sub-KZ-doctrine, thus with the components of $\varepsilon\:\T\tto\P$ being fully faithful.

A second remark: in the situation of Proposition \ref{a9.1}, the components of the transformation $\varepsilon\:\T\tto\P$ are necessarily given by pointwise left Kan extensions. More precisely, $\<Y_{\bbA},\eta_{\bbA}\>\:\T(\bbA)\to\P(\bbA)$ is the $\T(\bbA)(\eta_{\bbA}-,-)$-weighted colimit of $Y_{\bbA}$ (which exists because $\P(\bbA)$ is cocomplete), and can thus be computed as
$$\<Y_{\bbA},\eta_{\bbA}\>\:\T(\bbA)\to\P(\bbA)\:t\mapsto \T(\bbA)(\eta_{\bbA}-,t).$$
By fully faithfulness of $\varepsilon_{\bbA}\:\T(\bbA)\to\P(\bbA)$ and the Yoneda Lemma, we can compute that
$$\T(\bbA)(\eta_{\bbA}-,t)=\P(\bbA)(\varepsilon_{\bbA}\circ\eta_{\bbA}-,\varepsilon_{\bbA}(t))=\P(\bbA)(Y_{\bbA}-,\varepsilon_{\bbA}(t))=\varepsilon_{\bbA}(t).$$
Hence the component of $\varepsilon\:\T\tto\P$ at $\bbA\in\Cat(\Q)$ is necessarily the Kan extension $\<Y_{\bbA},\eta_{\bbA}\>$. We can push this argument a little further to obtain a characterisation of those KZ-doctrines which occur as full sub-KZ-doctrines of $\P$:
\begin{corollary}\label{a9.2}
A KZ-doctrine $(\T,\mu,\eta)$ on $\Cat(\Q)$ is a full sub-KZ-doctrine of $\P$ if and only if all $\eta_{\bbA}\:\bbA\to\T(\bbA)$ and all left Kan extensions $\<Y_{\bbA},\eta_{\bbA}\>\:\T(\bbA)\to\P(\bbA)$ are fully faithful. 
\end{corollary}
\proof
If $\T$ is a full sub-KZ-doctrine of $\P$, then we have just remarked that $\varepsilon_{\bbA}=\<Y_{\bbA},\eta_{\bbA}\>$, and thus these Kan extensions are fully faithful. Moreover -- because $\varepsilon_{\bbA}\circ\eta_{\bbA}=Y_{\bbA}$ with both $\varepsilon_{\bbA}$ and $Y_{\bbA}$ fully faithful -- also $\eta_{\bbA}$ must be fully faithful. 

Conversely, if $(\T,\mu,\eta)$ is a KZ-doctrine with each $\eta_{\bbA}\:\bbA\to\T(\bbA)$ fully faithful, then -- e.g.\ by [Stubbe, 2005, Proposition 6.7] -- the left Kan extensions $\<Y_{\bbA},\eta_{\bbA}\>$ (exist and) satisfy $\<Y_{\bbA},\eta_{\bbA}\>\circ\eta_{\bbA}\cong Y_{\bbA}$. By assumption each of these Kan extensions is fully faithful, so we must now prove that they are the components of a natural transformation and that this natural transformation commutes with the multiplications of $\T$ and $\P$. We do this in four steps:

(i) For any $\bbA\in\Cat(\Q)$, there is the free $\T$-algebra $\mu_{\bbA}\:\T(\T(\bbA))\to\T(\bbA)$. But the free $\P$-algebra $M_{\bbA}\:\P(\P(\bbA))\to\P(\bbA)$ on $\P(\bbA)$ also induces a $\T$-algebra on $\P(\bbA)$: namely, $ M_{\bbA}\circ\<Y_{\P(\bbA)},\eta_{\P(\bbA)}\>\:\T(\P(\bbA)\to\P(\bbA)$. To see this, it suffices to prove the adjunction $ M_{\bbA}\circ\<Y_{\P(\bbA)},\eta_{\P(\bbA)}\>\dashv\eta_{\P(\bbA)}$. The counit is easily checked:
$$ M_{\bbA}\circ\<Y_{\P(\bbA)},\eta_{\P(\bbA)}\>\circ\eta_{\P(\bbA)}= M_{\bbA}\circ Y_{\P(\bbA)}=1_{\P(\bbA)},$$
using first the factorisation property of the Kan extension and then the split adjunction $M_{\bbA}\dashv Y_{\P(\bbA)}$. As for the unit of the adjunction, we compute that
\begin{eqnarray*}
\eta_{\P(\bbA)}\circ M_{\bbA}\circ\<Y_{\P(\bbA)},\eta_{\P(\bbA)}\>
 & = & \T( M_{\bbA}\circ\<Y_{\P(\bbA)},\eta_{\P(\bbA)}\>)\circ\eta_{\T(\P(\bbA))} \\
 & \geq & \T( M_{\bbA}\circ\<Y_{\P(\bbA)},\eta_{\P(\bbA)}\>)\circ\T(\eta_{\P(\bbA)}) \\
 & = & \T( M_{\bbA}\circ\<Y_{\P(\bbA)},\eta_{\P(\bbA)}\>\circ\eta_{\P(\bbA)}) \\
 & = & \T(1_{\P(\bbA)}) \\
 & = & 1_{\T(\P(\bbA))},
\end{eqnarray*}
using naturality of $\eta$ and the KZ inequality for $\T$, and recycling the computation we made for the counit.

(ii) Next we prove, for each $\Q$-category $\bbA$, that $\<Y_{\bbA},\eta_{\bbA}\>\:\T(\bbA)\to\P(\bbA)$ is a $\T$-algebra homomorphism, for the algebra structures explained in the previous step. This is the case if and only if $\<Y_{\bbA},\eta_{\bbA}\>=( M_{\bbA}\circ\<Y_{\P(\bbA)},\eta_{\P(\bbA)}\>)\circ\T(\<Y_{\bbA},\eta_{\bbA}\>)\circ\eta_{\T(\bbA)}$ (because the domain of $\<Y_{\bbA},\eta_{\bbA}\>$ is a free $\T$-algebra), and indeed:
\begin{eqnarray*}
  & & M_{\bbA}\circ\<Y_{\P(\bbA)},\eta_{\P(\bbA)}\>\circ\T(\<Y_{\bbA},\eta_{\bbA}\>)\circ\eta_{\T(\bbA)} \\
 & & = M_{\bbA}\circ\<Y_{\P(\bbA)},\eta_{\P(\bbA)}\>\circ\eta_{\P(\bbA)}\circ\<Y_{\bbA},\eta_{\bbA}\> \\
 & & = M_{\bbA}\circ Y_{\P(\bbA)}\circ\<Y_{\bbA},\eta_{\bbA}\> \\
 & & = 1_{\P(\bbA)}\circ\<Y_{\bbA},\eta_{\bbA}\> \\
 & & = \<Y_{\bbA},\eta_{\bbA}\>.
\end{eqnarray*}

(iii) To check that the left Kan extensions are the components of a natural transformation we must verify, for any $F\:\bbA\to\bbB$ in $\Cat(\Q)$, that $\P(F)\circ\<Y_{\bbA},\eta_{\bbA}\>=\<Y_{\bbB},\eta_{\bbB}\>\circ\T(F)$. Since this is an equation of $\T$-algebra homomorphisms for the $\T$-algebra structures discussed in step (i) -- concerning $\P(F)$, it is easily seen to be a left adjoint and therefore also a $\T$-algebra homomorphism [Kock, 1995, Proposition 2.5] -- it suffices to show that $\P(F)\circ\<Y_{\bbA},\eta_{\bbA}\>\circ\eta_{\bbA}=\<Y_{\bbB},\eta_{\bbB}\>\circ\T(F)\circ\eta_{\bbA}$. This is straightforward from the factorisation property of the Kan extension and the naturality of $Y_{\bbA}$ and $\eta_{\bbA}$.

(iv) Finally, the very fact that $\<Y_{\bbA},\eta_{\bbA}\>\:\T(\bbA)\to\P(\bbA)$ is a $\T$-algebra homomorphism as in step (ii), means that 
$$\xymatrix@=18mm{
\T(\T(\bbA))\ar[r]^{\T(\<Y_{\bbA},\eta_{\bbA}\>)}\ar[d]_{\mu_{\bbA}} & \T(\P(\bbA))\ar[r]^{\<Y_{\P(\bbA)},\eta_{\P(\bbA)}\>} & \P(\P(\bbA))\ar[d]^{ M_{\bbA}} \\
\T(\bbA)\ar[rr]_{\<Y_{\bbA},\eta_{\bbA}\>} & & \P(\bbA)}$$
commutes: it expresses precisely the compatibility of the natural transformation whose components are the Kan extensions, with the multiplications of, respectively, $\T$ and $\P$.
\endofproof

\section{Interlude: classifying cotabulations}

In this section it is Proposition \ref{a8} which is of most interest: it explains in particular how the 2-functors on $\Cat(\Q)$ of Proposition \ref{a9.1} can be extended to $\Dist(\Q)$. It could easily be proved with a direct proof, but it seemed more appropriate to include first some material on classifying cotabulations, then use this to give a somewhat more conceptual proof of (the quantaloidal generalisation of) Akhvlediani {\it et al.}'s `Extension Theorem' [2009, Theorem 1] in our Proposition \ref{a7}, and finally derive Proposition \ref{a8} as a particular case.

A {\bf cotabulation} of a distributor $\Phi\:\bbA\dist\bbB$ between $\Q$-categories is a pair of functors, say $S\:\bbA\to\bbC$ and $T\:\bbB\to\bbC$, such that $\Phi=\bbC(T-,S-)$. If $F\:\bbC\to\bbC'$ is a fully faithful functor then also $F\circ S\:\bbA\to\bbC'$ and $F\circ T\:\bbB\to\bbC'$ cotabulate $\Phi$; so a distributor admits many different cotabulations. But the classifying property of $\P(\bbB)$ suggests a particular one:
\begin{proposition}\label{a4.1}
Any distributor $\Phi\:\bbA\dist\bbB$ is cotabulated by $Y_{\Phi}\:\bbA\to\P(\bbB)$ and $Y_{\bbB}\:\bbB\to\P(\bbB)$. We call this pair the {\bf classifying cotabulation} of $\Phi\:\bbA\dist\bbB$.
\end{proposition}
\proof
We compute for $a\in\bbA$ and $b\in\bbB$ that $\P(\bbB)(Y_{\bbB}(b),Y_{\Phi}(a))=Y_{\Phi}(a)(b)=\Phi(b,a)$ by using the Yoneda Lemma.
\endofproof

For two distributors $\Phi\:\bbA\dist\bbB$ and $\Psi\:\bbB\dist\bbC$ it is easily seen that $Y_{\Psi\tensor\Phi}=\P_0(\Psi)\circ Y_{\Phi}$, so the classifying cotabulation of the composite $\Psi\tensor\Phi$ relates to those of $\Phi$ and $\Psi$ as
\begin{equation}\label{a5}
\Psi\tensor\Phi=\P(\bbC)(\P_0(\Psi)\circ Y_{\Phi}-,Y_{\bbC}-).
\end{equation}
For a functor $F\:\bbA\to\bbB$ it is straightforward that $Y_{\bbB(-,F-)}=Y_{\bbB}\circ F$, so
\begin{equation}\label{a6}
\bbB(-,F-)=\P(\bbB)(Y_{\bbB}-,Y_{\bbB}\circ F-).
\end{equation}
In particular, the identity distributor $\bbA\:\bbA\dist\bbA$ has the classifying cotabulation 
\begin{equation}\label{a6.1}
\bbA=\P(\bbA)(Y_{\bbA}-,Y_{\bbA}-).
\end{equation}
Given that classifying cotabulations are thus perfectly capable of encoding composition and identities, it is natural to extend a given endo-functor on $\Cat(\Q)$ to an endo-functor on $\Dist(\Q)$ by applying it to classifying cotabulations. Now follows a statement of the `Extension Theorem' of [Akhvlediani {\it et al.}, 2009] in the generality of quantaloid-enriched category theory, and a proof based on the calculus of classifying cotabulations.
\begin{proposition}[Akhvlediani et al., 2009]\label{a7}
\ Any 2-functor $\T\:\Cat(\Q)\to\Cat(\Q)$ extends to a lax 2-functor $\T'\:\Dist(\Q)\to\Dist(\Q)$, which is defined to send a distributor $\Phi\:\bbA\dist\bbB$ to the distributor cotabulated by 
$\T(Y_{\Phi})\:\T(\bbA)\to\T(\P(\bbB))$ and $\T(Y_{\bbB})\:\T(\bbB)\to\T(\P(\bbB))$. This comes with a lax transformation
\begin{equation}\label{a7.1}
\begin{array}{c}
\xymatrix@=7ex{
\Dist(\Q)\ar[r]^{\T'} & \Dist(\Q) \\
\Cat(\Q)\ar[u]^i\ar[r]_{\T}\ar@{}[ur]|{\rotatebox{-45}{$\Longleftarrow$}} & \Cat(\Q)\ar[u]_i}
\end{array}
\end{equation}
all of whose components are identities. This lax transformation is a (strict) 2-natural transformation (i.e.\ this diagram is commutative) if and only if $\T'$ is normal, if and only if each $\T(Y_{\bbA})\:\T(\bbA)\to\T(\P(\bbA))$ is fully faithful.
\end{proposition}
\proof
If $\Phi\leq\Psi$ holds in $\Dist(\Q)(\bbA,\bbB)$ then (and only then) $Y_{\Phi}\leq Y_{\Psi}$ holds in $\Cat(\Q)(\bbA,\P(\bbB))$. By 2-functoriality of $\T\:\Cat(\Q)\to\Cat(\Q)$ we find that $\T(Y_{\Phi})\leq\T(Y_{\Psi})$, and thus $\T'(\Phi)\leq\T'(\Psi)$. 

Now suppose that $\Phi\:\bbA\dist\bbB$ and $\Psi\:\bbB\dist\bbC$ are given. Applying $\T$ to the commutative diagram
$$\xymatrix@=5ex{
\bbA\ar[dr]_{Y_{\Phi}} & & \bbB\ar[dl]_{Y_{\bbB}}\ar[dr]^{Y_{\Psi}} & & \bbC\ar[dl]_{Y_{\bbC}} \\
 & \P(\bbB)\ar[rr]_{\P_0(\Psi)} & & \P(\bbC)}$$
gives a commutative diagram in $\Cat(\Q)$, which embeds as a commutative diagram of left adjoints in the quantaloid $\Dist(\Q)$ by application of $i\:\Cat(\Q)\to\Dist(\Q)$. Lemma \ref{a2}, the formula in \eqref{a5} and the definition of $\T'$ allow us to conclude that $\T'(\Psi)\tensor\T'(\Phi)\leq\T'(\Psi\tensor\Phi)$.

Similarly, given $F\:\bbA\to\bbB$ in $\Cat(\Q)$, applying $\T$ to the commutative diagram
$$\xymatrix{
\bbA\ar[dr]_F & & \bbB\ar[dl]_{1_{\bbB}}\ar[dr]^{Y_{\bbB}} & & \bbB\ar[dl]^{Y_{\bbB}} \\
 & \bbB\ar[rr]_{Y_{\bbB}} & & \P(\bbB)}$$
gives a commutative diagram in $\Cat(\Q)$. This again embeds as a diagram of left adjoints in $\Dist(\Q)$ via $i\:\Cat(\Q)\to\Dist(\Q)$. Lemma \ref{a2}, the formula in \eqref{a6} and the definition of $\T'$ then straightforwardly imply that
\begin{eqnarray*}
\T'(\bbB(-,F-)) & = & \T(\P(\bbB))(\T(Y_{\bbB})-,\T(Y_{\bbB})-)\tensor\T(\bbB)(-,\T(F)-) \\
 & \geq & \T(\bbB)(-,\T(F)-),
\end{eqnarray*}
accounting for the lax transformation in \eqref{a7.1}. 

It further follows from this inequation, by applying it to identity functors, that $\T'$ is in general lax on identity distributors. But Lemma \ref{a2} also says: (i) if each $\T(Y_{\bbB})\:\T(\bbB)\to\T(\P(\bbB))$ is fully faithful (equivalently, if $\T'$ is normal), then necessarily $(i\circ\T)(F)\cong(\T'\circ i)(F)$ for all $F:\bbA\to\bbB$ in $\Cat(\Q)$, asserting that the diagram in \eqref{a7.1} commutes; (ii) and conversely, if that diagram commutes, then chasing the identities in $\Cat(\Q)$ shows that $\T'$ is normal.
\endofproof

We shall be interested in extending full sub-KZ-doctrines of the free cocompletion doctrine $\P\:\Cat(\Q)\to\Cat(\Q)$ to $\Dist(\Q)$; for this we make use of the functor $\P_1\:\Dist(\Q)\to\Dist(\Q)$ defined in the diagram in \eqref{a4}.
\begin{proposition}\label{a8}
Let $(\T,\varepsilon)$ be a full sub-KZ-doctrine of $\P\:\Cat(\Q)\to\Cat(\Q)$. 
The lax extension $\T'\:\Dist(\Q)\to\Dist(\Q)$ of $\T\:\Cat(\Q)\to\Cat(\Q)$ (as in Proposition \ref{a7}) can then be computed as follows: for $\Phi\:\bbA\dist\bbB$,
\begin{equation}\label{a8.1}
\T'(\Phi)=\P(\bbB)(\varepsilon_{\bbB}-,-)\tensor\P_1(\Phi)\tensor\P(\bbA)(-,\varepsilon_{\bbA}-).
\end{equation}
Moreover, $\T'$ is always a normal lax $\Sup$-functor, thus the diagram in \eqref{a7.1} commutes.
\end{proposition}
\proof
Let $\Phi\:\bbA\dist\bbB$ be a distributor. Proposition \ref{a7} defines $\T'(\Phi)$ to be the distributor cotabulated by $\T(Y_{\Phi})$ and $\T(Y_{\bbB})$; but by fully faithfulness of the components of $\varepsilon\:\T\tto\P$, and its naturality, we can compute that
\begin{eqnarray*}
 & & \T(\P(\bbB))(\T(Y_{\bbB})-,\T(Y_{\Phi})-) \\
 & & =\P(\P(\bbB))((\varepsilon_{\P(\bbB)}\circ\T(Y_{\bbB}))-,(\varepsilon_{\P(\bbB)}\circ\T(Y_{\Phi}))-) \\
 & & =\P(\P(\bbB))((\P(Y_{\bbB})\circ\varepsilon_{\bbB})-,(\P(Y_{\Phi})\circ\varepsilon_{\bbA})-) \\
 & & =\P(\bbB)(\varepsilon_{\bbB}-,-)\tensor\P(\P(\bbB))(\P(Y_{\bbB})-,\P(Y_{\Phi})-)\tensor\P(\bbA)(-,\varepsilon_{\bbA}-).
\end{eqnarray*}
The middle term in this last expression can be reduced:
\begin{eqnarray*}
\P(\P(\bbB))(\P(Y_{\bbB})-,\P(Y_{\Phi})-)
 & = & [\P(\bbB)(-,Y_{\bbB}-)\tensor-,\ \P(\bbB)(-,Y_{\Phi}-)\tensor-] \\
 & = & [-,\P(\bbB)(Y_{\bbB}-,-)\tensor\P(\bbB)(-,Y_{\Phi}-)\tensor-] \\
 & = & [-,\P(\bbB)(Y_{\bbB}-,Y_{\Phi}-)\tensor-] \\
 & = & [-,\Phi\tensor-] \\
 & = & \P(\bbB)(-,\P_0(\Phi)-) \\
 & = & (i\circ\P_0)(\Phi)(-,-) \\
 & = & \P_1(\Phi)(-,-).
\end{eqnarray*}
Thus we arrive at \eqref{a8.1}. Because $\P_1$ is a (strict) functor and because each $\varepsilon_{\bbA}$ is fully faithful, it follows from \eqref{a8.1} that $\T'$ is normal. Similarly, because $\P_1$ is a $\Sup$-functor, $\T'$ preserves local suprema too. 
\endofproof
If we apply Proposition \ref{a7} to the 2-functor $\P\:\Cat(\Q)\to\Cat(\Q)$ itself, then we find that $\P'=\P_1$ (and thus it is strictly functorial, not merely normal lax). In general however, $\T'$ does {\em not} preserve composition.

\section{Cocompletion: saturated classes of weights vs.\ KZ-doctrines} 

The $\Phi$-weighted colimit of a functor $F$ exists if and only if, for every $a\in\bbA_0$, $\colim(\Phi(-,a),F)$ exists:
$$\xymatrix@=15mm{
\bbB\ar[rr]^F & & \bbC \\
\bbA\ar[u]|{\distsign}_{\Phi}\ar@{.>}[rru]|{\colim(\Phi,F)} \\
{*_{ta}}\ar[u]|{\distsign}_{\bbA(-,a)}\ar@{.>}[uurr]_{\colim(\Phi(-,a),F)}\ar@<1mm>@/^5mm/[uu]|{\distsign}^{\Phi(-,a)=\Phi\tensor\bbA(-,a)} }$$
Indeed, $\colim(\Phi,F)(a)=\colim(\Phi(-,a),F)(*)$. But now $\Phi(-,a)\:*_{ta}\dist\bbB$ is a presheaf on $\bbB$. As a consequence, a $\Q$-category $\bbC$ is cocomplete if and only if it admits all colimits weighted by presheaves.

It therefore makes perfect sense to fix a class $\C$ of presheaves and study those $\Q$-categories that admit all colimits weighted by elements of $\C$: by definition these are the {\bf $\C$-cocomplete categories}.  Similarly, a functor $G\:\bbC\to\bbC'$ is {\bf $\C$-cocontinuous} if it preserves all colimits weighted by elements of $\C$.

As [Albert and Kelly, 1988; Kelly and Schmitt, 2005] demonstrated in the case of $\V$-categories (for $\V$ a symmetric monoidal closed category with locally small, complete and cocomplete underlying category $\V_0$), and as we shall argue here for $\Q$-categories too, it is convenient to work with classes of presheaves that ``behave nicely'':
\begin{definition}\label{a10}
A class $\C$ of presheaves on $\Q$-categories is {\bf saturated} if:
\begin{enumerate}
\item $\C$ contains all representable presheaves, 
\item for each $\phi\:*_X\dist\bbA$ in $\C$ and each functor $G\:\bbA\to\P(\bbB)$ for which each $G(a)$ is in $\C$, $\colim(\phi,G)$ is in $\C$ too.
\end{enumerate} 
\end{definition}
There is another way of putting this. Observe first that any class $\C$ of presheaves on $\Q$-categories defines a sub-2-graph $k\:\Dist_{\C}(\Q)\incl\Dist(\Q)$ by
\begin{equation}\label{a10.1}
\Phi\:\bbA\dist\bbB\mbox{ is in }\Dist_{\C}(\Q)\stackrel{{\rm def.}}{\Longleftrightarrow}\mbox{ for all $a\in\bbA_0$: }\Phi(-,a)\in\C.
\end{equation}
Then in fact we have:
\begin{proposition}\label{a11}
A class $\C$ of presheaves on $\Q$-categories is saturated if and only if $\Dist_{\C}(\Q)$ is a sub-2-category of $\Dist(\Q)$ containing (all objects and) all identities. In this case there is an obvious factorisation
$$\xymatrix{
\Cat(\Q)\ar[rr]^i\ar@{.>}[rd]_j & & \Dist(\Q) \\
 & \Dist_{\C}(\Q)\ar@{^{(}->}[ru]_k}$$
\end{proposition}
\proof
With \eqref{a10.1} it is trivial that $\C$ contains all representable presheaves if and only if $\Dist_{\C}(\Q)$ contains all objects and all identities. 

Next, assume that $\C$ is a saturated class of presheaves, and let $\Phi\:\bbA\dist\bbB$ and $\Psi\:\bbB\dist\bbC$ be arrows in $\Dist_{\C}(\Q)$. Invoking the classifying property of $\P(\bbC)$ and the computation of colimits in $\P(\bbC)$, we find $\colim(\Phi(-,a),Y_{\Psi})=\Psi\tensor\Phi(-,a)$ for each $a\in\bbA_0$. But because $\Phi(-,a)\in\C$ and for each $b\in\bbB_0$ also $Y_{\Psi}(b)=\Psi(-,b)\in\C$, this colimit, i.e.\ $\Psi\tensor\Phi(-,a)$, is an element of $\C$. This holds for all $a\in\bbA_0$, thus the composition $\Psi\tensor\Phi\:\bbA\dist\bbC$ is an arrow in $\Dist_{\C}(\Q)$.

Conversely, assuming $\Dist_{\C}(\Q)$ is a sub-2-category of $\Dist(\Q)$, let $\phi\:*_A\dist\bbB$ be in $\C$ and let $F\:\bbB\to\P(\bbC)$ be a functor such that, for each $b\in\bbB$, $F(b)$ is in $\C$. By the classifying property of $\P(\bbC)$ we can equate the functor $F\:\bbB\to\P(\bbC)$ with a distributor $\Phi_F\:\bbB\dist\bbC$ and by the computation of colimits in $\P(\bbC)$ we know that $\colim(\phi,F)=\Phi_F\tensor\phi$. Now $\Phi_F(-,b)=F(b)$ by definition, so $\Phi_F\:\bbB\dist\bbC$ is in $\Dist_{\C}(\Q)$; but also $\phi:*_A\dist\bbB$ is in $\Dist_{\C}(\Q)$, and therefore their composite is in $\Dist_{\C}(\Q)$, i.e.\ $\colim(\phi,F)$ is in $\C$, as wanted.

Finally, if $F\:\bbA\to\bbB$ is any functor, then for each $a\in\bbA$ the representable $\bbB(-,Fa)\:*_{ta}\dist\bbB$ is in the saturated class $\C$, and therefore $\bbB(-,F-)\:\bbA\dist\bbB$ is in $\Dist_{\C}(\Q)$. This accounts for the factorisation of $\Cat(\Q)\to\Dist(\Q)$ over $\Dist_{\C}(\Q)\incl\Dist(\Q)$.
\endofproof

We shall now characterise saturated classes of presheaves on $\Q$-categories in terms of KZ-doctrines on $\Cat(\Q)$. (We shall indeed always deal with a saturated class of presheaves, even though certain results hold under weaker hypotheses.) We begin by pointing out a classifying property:
\begin{proposition}\label{a13}
Let $\C$ be a saturated class of presheaves and, for a $\Q$-category $\bbA$, write $J_{\bbA}\:\C(\bbA)\to\P(\bbA)$ for the full subcategory of $\P(\bbA)$ determined by those presheaves on $\bbA$ which are elements of $\C$. A distributor $\Phi\:\bbA\dist\bbB$ belongs to $\Dist_{\C}(\Q)$ if and only if there exists a (necessarily unique) factorisation
\begin{equation}\label{a13.1}
\begin{array}{c}
\xymatrix{
\bbA\ar[rr]^{Y_{\Phi}}\ar@{.>}[dr]_{I_{\Phi}} & & \P(\bbB) \\
 & \C(\bbB)\ar[ru]_{J_{\bbB}}}
\end{array}
\end{equation}
in which case $\Phi$ is cotabulated by $I_{\Phi}\:\bbA\to\C(\bbB)$ and $I_{\bbB}\:\bbB\to\C(\bbB)$ (the latter being the factorisation of $Y_{\bbB}$ through $J_{\bbB}$).
\end{proposition}
\proof
The factorisation property in \eqref{a13.1} literally says that, for any $a\in\bbA$, the presheaf $Y_{\Phi}(a)$ on $\bbB$ must be an element of the class $\C$. But $Y_{\Phi}(b)=\Phi(-,b)$ hence this is trivially equivalent to the statement in \eqref{a10.1}, defining those distributors that belong to $\Dist_{\C}(\Q)$. In particular, if $\C$ is saturated then $\Dist_{\C}(\Q)$ contains all identities, hence we have factorisations $Y_{\bbB}=J_{\bbB}\circ I_{\bbB}$ of the Yoneda embeddings. Hence, whenever a factorisation as in \eqref{a13.1} exists, we can use the fully faithful $J_{\bbB}\:\C(\bbB)\to\P(\bbB)$ to compute, starting from the classifying cotabulation of $\Phi$, that
$$\Phi=\P(\bbB)(Y_{\bbB}-,Y_{\Phi}-)=\P(\bbB)(J_{\bbB}(I_{\bbB}(-)),J_{\bbB}(I_{\Phi}(-)))=\C(\bbB)(I_{\bbB}-,I_{\Phi}-),$$
confirming the cotabulation of $\Phi$ by $I_{\Phi}$ and $I_{\bbB}$.
\endofproof

Any saturated class $\C$ thus automatically comes with the 2-functor
$$\C_0\:\Dist_{\C}(\Q)\to\Cat(\Q)\:\Big(\Phi\:\bbA\dist\bbB\Big)\mapsto\Big(\Phi\tensor-\:\C(\bbA)\to\C(\bbB)\Big)$$
and the full embeddings $J_{\bbA}\:\C(\bbA)\to\P(\bbA)$ are the components of a 2-natural transformation
$$\xy\xymatrix@R=5mm{
 & \Dist(\Q)\ar@/^2mm/[rd]^{\P_0} \\
\Dist_{\C}(\Q)\ar@/^2mm/[ur]^k\ar[rr]_{\C_0} & & \Cat(\Q)}
\POS(23,-6)\drop{\rotatebox{90}{$\Longrightarrow$}}
\POS(26,-6)\drop{J}
\endxy$$
Composing $\C_0$ with $j\:\Cat(\Q)\to\Dist_{\C}(\Q)$ it is natural to define
$$\C\:\Cat(\Q)\to\Cat(\Q)\:\Big(F\:\bbA\to\bbB\Big)\mapsto\Big(\bbB(-,F-)\tensor-\:\C(\bbA)\to\C(\bbB)\Big)$$
together with
$$\xy\xymatrix@R=5mm{
\Cat(\Q)\ar@/^4mm/[r]^{\P}\ar@/_4mm/[r]_{\C} & \Cat(\Q)}
\POS(11,0)\drop{\rotatebox{90}{$\Longrightarrow$}}
\POS(13,0)\drop{J}
\endxy$$
(slightly abusing notation). We apply previous results, particularly Proposition \ref{a9.1}:
\begin{proposition}\label{a16}
If $\C$ is a saturated class of presheaves on $\Q$-categories then the 2-functor $\C\:\Cat(\Q)\to\Cat(\Q)$ together with the transformation $J\:\C\tto\P$ forms a full sub-KZ-doctrine of $\P$. Moreover, the $\C$-cocomplete $\Q$-categories are precisely the $\C$-algebras, and the $\C$-cocontinuous functors between $\C$-cocomplete $\Q$-categories are precisely the $\C$-algebra homomorphisms.
\end{proposition}
\proof 
To fulfill the hypotheses in Proposition \ref{a9.1}, we only need to check the factorisation of the multiplication: if we prove, for any $\Q$-category $\bbA$ and each $\phi\in\C(\C(\bbA))$, that the $(J*J)_{\bbA}(\phi)$-weighted colimit of $1_{\P(\bbA)}$ is in $\C(\bbA)$, then we obtain the required commutative diagram
$$\xymatrix@C=15ex@R=10mm{
\P(\P(\bbA))\ar[r]^{\colim(-,1_{\P(\bbA)})} & \P(\bbA) \\
\C(\C(\bbA))\ar@{.>}[r]\ar[u]^{(J*J)_{\bbA}} & \C(\bbA)\ar[u]_{J_{\bbA}}}$$
But because $(J*J)_{\bbA}=\P(J_{\bbA})\circ J_{\C(\bbA)}$ we can compute that
$$\colim((J*J)_{\bbA}(\phi),1_{\P(\bbA)})=\colim(\P(\bbA)(-,J_{\bbA}-)\tensor\phi,1_{\P(\bbA)})=\colim(\phi,J_{\bbA})$$
and this colimit indeed belongs to the saturated class $\C$, because both $\phi$ and (the objects in) the image of $J_{\bbA}$ are in $\C$.

A $\Q$-category $\bbB$ is a $\C$-algebra if and only if $I_{\bbB}\:\bbB\to\C(\bbB)$ admits a left adjoint in $\Cat(\Q)$ (because $\C$ is a KZ-doctrine). Suppose that $\bbB$ is indeed a $\C$-algebra, and write the left adjoint as $L_{\bbB}\:\C(\bbB)\to\bbB$. If $\phi\:*_X\dist\bbA$ is a presheaf in $\C$ and $F\:\bbA\to\bbB$ is any functor, then $\C(F)(\phi)$ is an object of $\C(\bbB)$, thus we can consider the object $L_{\bbB}(\C(F)(\phi))$ of $\bbB$. This is precisely the $\phi$-weighted colimit of $F$, for indeed its universal property holds: for any $b\in\bbB$,
\begin{eqnarray*}
\bbB(L_{\bbB}(\C(F)(\phi)),b) 
 & = & \C(\bbB)(\C(F)(\phi),I_{\bbB}(b)) \\
 & = & \P(\bbB)(J_{\bbB}(\C(F)(\phi)),J_{\bbB}(I_{\bbB}(b))) \\
 & = & [\P(F)(J_{\bbB}(\phi)),Y_{\bbB}(b)] \\
 & = & [\bbB(-,F-)\tensor J_{\bbB}(\phi),\bbB(-,b)] \\
 & = & [J_{\bbB}(\phi),\bbB(F-,-)\tensor\bbB(-,b)] \\
 & = & [\phi,\bbB(F-,b)].
\end{eqnarray*}
(Apart from the adjunction $L_{\bbB}\dashv I_{\bbB}$ we used the fully faithfulness of $J_{\bbB}$ and its naturality, and then made some computations with liftings and adjoints in $\Dist(\Q)$.)

Conversely, suppose that $\bbB$ admits all $\C$-weighted colimits. In particular can we then compute, for any $\phi\in\C(\bbB)$, the $\phi$-weighted colimit of $1_{\bbB}$, and doing so gives a function $f\:\C(\bbB)\to\bbB\:\phi\mapsto\colim(\phi,1_{\bbB})$. But for any $\phi\in\C(\bbB)$ and any $b\in\bbB$ it is easy to compute, from the universal property of colimits and using the fully faithfulness of $J_{\bbB}$, that
\begin{eqnarray*}
\bbB(f(\phi),b)
 & = & [\phi,\bbB(1_{\bbB}-,b)] \ = \ \P(\bbB)(\phi,Y_{\bbB}(b)) \\
 & = & \P(\bbB)(J_{\bbB}(\phi),J_{\bbB}(I_{\bbB}(b))) \ = \ \C(\bbB)(\phi,I_{\bbB}(b)).
\end{eqnarray*}
This straightforwardly implies that $\phi\mapsto f(\phi)$ is in fact a functor (and not merely a function), and that it is left adjoint to $I_{\bbB}$; thus $\bbB$ is a $\C$-algebra.

Finally, let $G\:\bbB\to\bbC$ be a functor between $\C$-cocomplete $\Q$-categories. Supposing that $G$ is $\C$-cocontinuous, we can compute any $\psi\in\C(\bbB)$ that
$$G(L_{\bbB}(\psi))=G(\colim(\psi),1_{\bbB})=\colim(\psi),G)=L_{\bbC}(\C(G)(\psi)),$$
proving that $G$ is a homomorphism between the $\C$-algebras $(\bbB,L_{\bbB})$ and $(\bbC,L_{\bbC})$. Conversely, supposing now that $G$ is a homomorphism, we can compute for any presheaf $\phi\:*_X\dist\bbA$ in $\C$ and any functor $F\:\bbA\to\bbB$ that
$$G(\colim(\phi,F))=G(L_{\bbB}(\C(F)(\phi)))=L_{\bbC}(\C(G)(\C(F)(\phi)))=\colim(\phi,G\circ F),$$
proving that $G$ is $\C$-cocontinuous.
\endofproof
Also the converse of the previous Proposition is true:
\begin{proposition}\label{a16.1}
If $(\T,\varepsilon)$ is a full sub-KZ-doctrine of $\P\:\Cat(\Q)\to\Cat(\Q)$ then
\begin{equation}\label{a16.2}
\C_{\T}:=\{\varepsilon_{\bbA}(t)\mid \bbA\in\Cat(\Q),\ t\in\T(\bbA)\}
\end{equation}
is a saturated class of presheaves on $\Q$-categories. Moreover, the $\T$-algebras are precisely the $\C_{\T}$-cocomplete categories, and the $\T$-algebra homomorphisms are precisely the $\C_{\T}$-cocontinuous functors between the $\C_{\T}$-cocomplete categories.
\end{proposition}
\proof
We shall write $\Dist_{\T}(\Q)$ for the sub-2-graph of $\Dist(\Q)$ determined -- as prescribed in \eqref{a10.1} -- by the class $\C_{\T}$, and we shall show that it is a sub-2-category containing all (objects and) identities of $\Dist(\Q)$. But a distributor $\Phi\:\bbA\dist\bbB$ belongs to $\Dist_{\T}(\Q)$ if and only if the classifying functor $Y_{\Phi}\:\bbA\to\P(\bbB)$ factors (necessarily essentially uniquely) through the fully faithful $\varepsilon_{\bbB}\:\T(\bbB)\to\P(\bbB)$. 

By hypothesis there is a factorisation $Y_{\bbA}=\varepsilon_{\bbA}\circ\eta_{\bbA}$ for any $\bbA\in\Cat(\Q)$, so $\Dist_{\T}(\Q)$ contains all identities. Secondly, suppose that $\Phi\:\bbA\dist\bbB$ and $\Psi\:\bbB\dist\bbC$ are in $\Dist_{\T}(\Q)$, meaning that there exist factorisations
$$\xymatrix@=6mm{
\bbA\ar[rr]^{Y_{\Phi}}\ar@{.>}[dr]_{I_{\Phi}} & & \P(\bbB) \\
 & \T(\bbB)\ar[ur]_{\varepsilon_{\bbB}}}
\hspace{1cm}
\xymatrix@=6mm{
\bbB\ar[rr]^{Y_{\Psi}}\ar@{.>}[dr]_{I_{\Psi}} & & \P(\bbC) \\
 & \T(\bbC)\ar[ur]_{\varepsilon_{\bbC}}}$$
The following diagram is then easily seen to commute:
$$\xymatrix@=15mm{
 & \T(\bbB)\ar[r]^{\varepsilon_{\bbB}}\ar[d]_{\T(I_{\Psi})} & \P(\bbB)\ar[d]^{\P(I_{\Psi})}\ar@<4mm>@/^15mm/[dd]^{\P(Y_{\Psi})} \\
 & \T(\T(\bbC))\ar[r]^{\varepsilon_{\T(\bbC)}}\ar[d]_{\T(\varepsilon_{\bbC})}\ar@{.>}[dr]|{(\varepsilon{*}\varepsilon)_{\bbC}}\ar[dl]_{\mu_{\bbC}} & \P(\T(\bbC))\ar[d]^{\P(\varepsilon_{\bbC})} \\
\T(\bbC)\ar[dr]_{\varepsilon_{\bbC}} & \T(\P(\bbC))\ar[r]_{\varepsilon_{\P(\bbC)}} & \P(\P(\bbC))\ar[dl]^{ M_{\bbC}} \\
 & \P(\bbC) }$$
But we can compute, for any $\phi\in\P(\bbB)$, that
\begin{eqnarray*}
( M_{\bbC}\circ\P(Y_{\Psi}))(\phi)
 & = & \colim(\P(\bbA)(-,Y_{\Psi}-)\tensor\phi,1_{\P(\bbA)}) \\
 & = & \colim(\phi,Y_{\Psi}) \\
 & = & \Psi\tensor\phi \\
 & = & \P_0(\Psi)(\phi)
\end{eqnarray*}
and therefore $Y_{\Psi\tensor\Psi}=\P_0(\Psi)\circ Y_{\Phi}= M_{\bbC}\circ\P(Y_{\Psi})\circ\varepsilon_{\bbB}\circ I_{\Phi}=\varepsilon_{\bbC}\circ\mu_{\bbC}\circ\T(I_{\Psi})\circ I_{\Phi}$, giving a factorisation of $Y_{\Psi\tensor\Phi}$ through $\varepsilon_{\bbC}$, as wanted.

The arguments to prove that a $\Q$-category $\bbB$ is a $\T$-algebra if and only if it is $\C_{\T}$-cocomplete, and that a $\T$-algebra homomorphism is precisely a $\C_{\T}$-cocontinuous functor between $\C_{\T}$-co\-comp\-lete $\Q$-categories, are much like those in the proof of Proposition \ref{a16}. Omitting the calculations, let us just indicate that for a $\T$-algebra $\bbB$, thus with a left adjoint $L_{\bbB}\:\T(\bbB)\to\bbB$ to $\eta_{\bbB}$, for any weight $\phi\:*_X\dist\bbA$ in $\C_{\T}$ -- i.e.\ $\phi=\varepsilon_{\bbA}(t)$ for some $t\in\T(\bbA)$ -- and any functor $F\:\bbA\to\bbB$, the object $L_{\bbB}(\T(F))(t)$ is the $\phi$-weighted colimit of $F$. And conversely, if $\bbB$ is a $\C_{\T}$-cocomplete $\Q$-category, then $\T(\bbB)\to\bbB\:t\mapsto\colim(\varepsilon_{\bbB}(t),1_{\bbB})$ is the left adjoint to $\eta_{\bbB}$, making $\bbB$ a $\T$-algebra.
\endofproof

If $\C$ is a saturated class of presheaves and we apply Proposition \ref{a16} to obtain a full sub-KZ-doctrine $(\C,J)$ of $\P\:\Cat(\Q)\to\Cat(\Q)$, then the application of Proposition \ref{a16.1} gives us back precisely that same class $\C$ that we started from. The other way round is slightly more subtle: if $(\T,\varepsilon)$ is a full sub-KZ-doctrine of $\P$ then Proposition \ref{a16.1} gives us a saturated class $\C_{\T}$ of presheaves, and this class in turn determines by Proposition \ref{a16} a full KZ-doctrine of $\P$, let us write it as $(\T',\varepsilon')$, which is {\em equivalent} to $\T$. More exactly, each (fully faithful) $\varepsilon_{\bbA}\:\T(\bbA)\to\P(\bbA)$ factors over the fully faithful and injective $\varepsilon'_{\bbA}\:\T'(\bbA)\to\P(\bbA)$, and this factorisation is fully faithful and surjective, thus an equivalence. These equivalences are the components of a 2-natural transformation $\delta\:\T\tto\T'$ which commutes with $\varepsilon$ and $\varepsilon'$.

We summarise all the above in the following:
\begin{theorem}\label{a16.3}
Propositions \ref{a16} and \ref{a16.1} determine an essentially bijective correspondence between, on the one hand, saturated classes $\C$ of presheaves on $\Q$-categories, and on the other hand, full sub-KZ-doctrines $(\T,\varepsilon)$ of the free cocompletion KZ-doctrine $\P\:\Cat(\Q)\to\Cat(\Q)$; a class $\C$ and a doctrine $\T$ correspond with each other if and only if the $\T$-algebras and their homomorphisms are precisely the $\C$-cocomplete $\Q$-categories and the $\C$-cocontinuous functors between them. Proposition \ref{a8} implies that, in this case, there is a normal lax $\Sup$-functor $\T'\:\Dist(\Q)\to\Dist(\Q)$, sending a distributor $\Phi\:\bbA\dist\bbB$ to the distributor $\T'(\Phi)\:\T(\bbA)\dist\T(\bbB)$ with elements
$$\T'(\Phi)(t,s)=\P(\bbB)(\varepsilon_{\bbB}(t),\Phi\tensor\varepsilon_{\bbA}(s))\mbox{, \ \ for $s\in\T(\bbA)$, $t\in\T(\bbB)$,}$$
which makes the following diagram commute:
$$\xymatrix{
\Dist(\Q)\ar[r]^{\T'} & \Dist(\Q) \\
\Cat(\Q)\ar[r]_{\T}\ar[u]^i & \Cat(\Q)\ar[u]_i}$$
\end{theorem}

\section{Conical cocompletion and the Hausdorff doctrine} 

\subsection{Conical colimits}

Let $\bbA$ be a $\Q$-category. Putting, for any $a,a'\in\bbA$,
$$a\leq a'\stackrel{\mbox{\footnotesize def.}}{\iff}ta=ta'\mbox{ and }1_{ta}\leq\bbA(a,a')$$
defines an order relation on the objects of $\bbA$. (There are equivalent conditions in terms of representable presheaves.) For a given $\Q$-category $\bbA$ and a given object $X\in\Q_0$, we shall write $(\bbA_X,\leq_X)$ for the ordered set of objects of $\bbA$ of type $X$. Because elements of different type in $\bbA$ can never have a supremum in $(\bbA_0,\leq)$, it would be very restrictive to require this order to admit arbitrary suprema; instead, experience shows that it makes good sense to require each $(\bbA_X,\leq_X)$ to be a sup-lattice: we then say that $\bbA$ is {\bf order-cocomplete} [Stubbe, 2006]. As spelled out in that reference, we have:
\begin{proposition}\label{2}
For a family $(a_i)_{i\in I}$ in $\bbA_X$, the following are equivalent:
\begin{enumerate}
\item $\bigvee_ia_i$ exists in $\bbA_X$ and $\bbA(\bigvee_ia_i,-)=\bigwedge_i\bbA(a_i,-)$ holds in $\Dist(\Q)(\bbA,*_X)$,
\item $\bigvee_ia_i$ exists in $\bbA_X$ and $\bbA(-,\bigvee_ia_i)=\bigvee_i\bbA(-,a_i)$ holds in $\Dist(\Q)(*_X,\bbA)$,
\item if we write $(I,\leq)$ for the ordered set in which $i\leq j$ precisely when $a_i\leq_X a_j$ and $\bbI$ for the free $\Q(X,X)$-category on the poset $(I,\leq)$, $F\:\bbI\to\bbA$ for the functor $i\mapsto a_i$ and $\gamma\:*_X\dist\bbI$ for the presheaf with values $\gamma(i)=1_X$ for all $i\in\bbI$, then the $\gamma$-weighted colimit of $F$ exists.
\end{enumerate}
In this case, $\colim(\gamma,F)=\bigvee_ia_i$ and it is the {\bf conical colimit} of $(a_i)_{i\in I}$ in $\bbA$.
\end{proposition} 
It is important to realise that such conical colimits -- which are enriched colimits!\ -- can be characterised by a property of weights:
\begin{proposition}\label{3}
For a presheaf $\phi\:*_X\dist\bbA$, the following conditions are equivalent:
\begin{enumerate}
\item there exists a family $(a_i)_{i\in I}$ in $\bbA_X$ such that for any functor $G\:\bbA\to\bbB$, if the $\phi$-weighted colimit of $G$ exists, then it is the conical colimit of the family $(G(a_i))_i$,
\item there exists a family $(a_i)_{i\in I}$ in $\bbA_X$ for which $\phi=\bigvee_i\bbA(-,a_i)$ holds in $\Dist(\Q)(*_X,\bbA)$,
\item there exist an ordered set $(I,\leq)$ and a functor $F\:\bbI\to\bbA$ with domain the free $\Q(X,X)$-category on $(I,\leq)$ such that, if we write $\gamma\:*_X\dist\bbI$ for the presheaf with values $\gamma(i)=1_X$ for all $i\in\bbI$, then $\phi=\bbA(-,F-)\tensor\gamma$.
\end{enumerate}
In this case, we call $\phi$ a {\bf conical presheaf}.
\end{proposition}
\proof
(i$\Rightarrow$ii) Applying the hypothesis to the functor $Y_{\bbA}\:\bbA\to\P(\bbA)$ -- indeed $\colim(\phi,Y_{\bbA})$ exists, and is equal to $\phi$ by the Yoneda Lemma -- we find a family $(a_i)_{i\in I}$ such that $\phi$ is the conical colimit in $\P(\bbA)$ of the family $(Y_{\bbA}(a))_i$. This implies in particular that $\phi=\bigvee_i\bbA(-,a_i)$.

(ii$\Rightarrow$iii) For $\phi=\bigvee_i\bbA(-,a_i)$ it is always the case that $[\bigvee_i\bbA(-,a_i),-]=\bigwedge_i[\bbA(-,a_i),-]$, i.e.\ $\P(\bbA)(\bigvee_i\bbA(-,a_i),-)=\bigwedge_i\P(\bbA)(\bbA(-,a_i),-)$. Thus $\phi$ is the conical colimit in $\P(\bbA)$ of the family $(\bbA(-,a_i))_i$, and Proposition \ref{2} allows for the conclusion.

(iii$\Rightarrow$i) If, for some functor $G\:\bbA\to\bbB$, $\colim(\phi,G)$ exists, then, by the hypothesis that $\phi=\bbA(-,F-)\tensor\gamma$, it is equal to $\colim(\bbA(-,F-)\tensor\gamma,G)=\colim(\gamma,G\circ F)$. The latter is the conical colimit of the family $(G(F(i)))_{i\in I}$; thus the family $(F(i))_i$ fulfills the requirement.
\endofproof
A warning is in order. Proposition \ref{3} attests that the conical presheaves on a $\Q$-category $\bbA$ are those which are a supremum of some family of representable presheaves on $\bbA$. Of course, neither that family of representables, nor the family of representing objects in $\bbA$, need to be unique.

Now comes the most important observation concerning conical presheaves.
\begin{proposition}\label{4}
The class of conical presheaves is saturated.
\end{proposition}
\proof
We shall check both conditions in Proposition \ref{a10}. 
All representable pre\-sheaves are clearly conical, so the first condition is fulfilled. As for the second condition, consider a conical presheaf $\phi\:*_X\dist\bbA$ and a functor $G\:\bbA\to\P(\bbB)$ such that each $G(a)\:*_{ta}\dist\bbB$ is a conical presheaf too. The $\phi$-weighted colimit of $G$ certainly exists, hence the first statement in Proposition \ref{3} applies: it says that $\colim(\phi,G)$ is the conical colimit of a family of conical presheaves. In other words, $\colim(\phi,G)$ is a supremum of a family of suprema of representables, and is therefore a supremum of representables too, hence a conical presheaf.
\endofproof

\subsection{The Hausdorff doctrine}

Applying Theorem \ref{a16.3} to the class of conical presheaves we get: 
\begin{definition}\label{5}
We write $\H\:\Cat(\Q)\to\Cat(\Q)$ for the KZ-doctrine associated with the class of conical presheaves. We call it the {\bf Hausdorff docrine on $\Cat(\Q)$}, and we say that $\H(\bbA)$ is the {\bf Hausdorff $\Q$-category} associated to a $\Q$-category $\bbA$. We write $\H'\:\Dist(\Q)\to\Dist(\Q)$ for the normal lax $\Sup$-functor which extends $\H$ from $\Cat(\Q)$ to $\Dist(\Q)$.
\end{definition}

To justify this terminology, and underline the concordance with [Akhvlediani {\it et al.}, 2009], we shall make this more explicit. According to Proposition \ref{a16}, $\H(\bbA)$ is the full subcategory of $\P(\bbA)$ determined by the conical presheaves on $\bbA$. By Proposition \ref{3} however, the objects of $\H(\bbA)$ can be equated with suprema of representables; so suppose that 
$$\phi=\bigvee_{a\in A}\bbA(-,a)\mbox{ \ \ and \ \ }\phi'=\bigvee_{a'\in A'}\bbA(-,a')$$
for subsets $A\subseteq\bbA_X$ and $A'\subseteq\bbA_Y$. Then we can compute that
\begin{eqnarray*}
\H(\bbA)(\phi',\phi)
 & = & \P(\bbA)(\phi',\phi) \\
 & = & [\phi',\phi] \\
 & = & [\bigvee_{a'}\bbA(-,a'),\bigvee_a\bbA(-,a] \\
 & = & \bigwedge_{a'}[\bbA(-,a'),\bigvee_a\bbA(-,a)] \\
 & = & \bigwedge_{a'}\bigvee_a[\bbA(-,a'),\bbA(-,a)] \\
 & = & \bigwedge_{a'}\bigvee_a\bbA(a',a).
\end{eqnarray*}
(The penultimate equality is due to the fact that each $\bbA(-,a')\:*_Y\dist\bbA$ is a left adjoint in the quantaloid $\Dist(\Q)$, and the last equality is due to the Yoneda lemma.) This is precisely the expected formula for the ``Hausdorff distance between (the conical presheaves determined by) the subsets $A$ and $A'$ of $\bbA$''. It must be noted that [Schmitt, 2006, Proposition 3.42] describes a very similar situation particularly for symmetric categories enriched in the commutative quantale of postive real numbers.

Similarly for functors: given a functor $F\:\bbA\to\bbB$ between $\Q$-categories, the functor $\H(F)\:\H(\bbA)\to\H(\bbB)$ sends a conical presheaf $\phi$ on $\bbA$ to the conical pre\-sheaf $\bbB(-,F-)\tensor\phi$ on $\bbB$. Supposing that $\phi=\bigvee_{a\in A}\bbA(-,a)$ for some $A\subseteq\bbA_X$, it is straightforward to check that
\begin{eqnarray*}
\bbB(-,F-)\tensor\phi
 & = & \bigvee_{x\in\bbA}\left(\bbB(-,Fx)\circ\bigvee_{a\in A}\bbA(x,a)\right) \\
 & = & \bigvee_{a\in A}\left(\bigvee_{x\in\bbA}\bbB(-,Fx)\circ\bbA(x,a)\right) \\
 & = & \bigvee_{a\in A}\bbB(-,Fa).
 \end{eqnarray*}
That is to say, ``$\H(F)$ sends (the conical presheaf determined by) $A\subseteq\bbA$ to (the conical presheaf determined by) $F(A)\subseteq\bbB$''.

Finally, by Proposition \ref{a8}, the action of $\H'$ on a distributor $\Phi\:\bbA\dist\bbB$ gives a distributor $\H'(\Phi)\:\H(\bbA)\dist\H'(\bbB)$ whose value in $\phi\in\H(\bbA)$ and $\psi\in\H(\bbB)$ is $\P(\bbB)(\psi,\Phi\tensor\phi)$. Assuming that
$$\phi=\bigvee_{a\in A}\bbA(-,a)\mbox{ and }\psi=\bigvee_{b\in B}\bbB(-,b)$$ 
for some $A\subseteq\bbA_X$ and $B\subseteq\bbB_Y$, a similar computation as above shows that 
$$\H'(\Phi)(\psi,\phi)=\bigwedge_{b\in B}\bigvee_{a\in A}\Phi(b,a).$$ 
This is the expected generalisation of the previous formula, to measure the ``Hausdorff distance between (the conical presheaves determined by) $A\subseteq\bbA$ and $B\subseteq\bbB$ through $\Phi\:\bbA\dist\bbB$''.

\subsection{Other examples}

The following examples of saturated classes of presheaves have been considered by [Kelly and Schmitt, 2005] in the case of categories enriched in symmetric monoidal categories.

\begin{example}[Minimal and maximal class]
The smallest saturated class of pre\-sheaves on $\Q$-categories is, of course, that containing only representable presheaves. It is straightforward that the KZ-doctrine on $\Cat(\Q)$ corresponding with this class is the identity functor. On the other hand, the class of all presheaves on $\Q$-categories corresponds with the free cocompletion KZ-doctrine on $\Cat(\Q)$.
\end{example}

\begin{example}[Cauchy completion]
The class of all left adjoint presheaves, also known as {\bf Cauchy presheaves}, on $\Q$-categories is saturated. Indeed, all representable presheaves are left adjoints. And suppose that $\Phi\:\bbA\dist\bbB$ and $\Psi\:\bbB\dist\bbC$ are distributors such that, for all $a\in\bbA$ and all $b\in\bbB$, $\Phi(-,a)\:*_{ta}\dist\bbB$ and $\Psi(-,b)\:*_{tb}\dist\bbC$ are left adjoints. Writing $\rho_b\:\bbC\dist*_{tb}$ for the right adjoint to $\Psi(-,b)$, it is easily verified that $\Psi$ is left adjoint to $\bigvee_{b\in\bbB}\bbB(-,b)\tensor\rho_b$. This makes sure that $(\Psi\tensor\Phi)(-,a)=\Psi\tensor\Phi(-,a)$ is a left adjoint too, and by Proposition \ref{a11} we can conclude that the class of Cauchy presheaves is saturated. The KZ-doctrine on $\Cat(\Q)$ which corresponds to this saturated class of presheaves, sends a $\Q$-category $\bbA$ to its {\bf Cauchy completion} [Lawvere, 1973; Walters, 1981; Street, 1983]. 
\end{example}

Inspired by the examples in [Lawvere, 1973] and the general theory in [Kelly and Schmitt, 2005], Vincent Schmitt [2006] has studied several other classes of presheaves for ordered sets (viewed as categories enriched in the 2-element Boolean algebra) and for generalised metric spaces (viewed as categories enriched in the quantale of positive real numbers). He constructs {\em saturated} classes of presheaves by requiring that each element of the class ``commutes'' (in a suitable way) with all elements of a given ({\em not-necessarily saturated}) class of presheaves. These interesting examples  do not seem to generalise straightforwardly to general quantaloid-enriched categories, so we shall not survey them here, but refer instead to [Schmitt, 2006] for more details.


\vspace{5mm}
\noindent
Isar Stubbe \\
Laboratoire de Math\'ematiques Pures et Appliqu\'ees \\
Universit\'e du Littoral-C\^ote d'Opale \\
50, rue Ferdinand Buisson \\
62228 Calais (France) \\
isar.stubbe@lmpa.univ-littoral.fr

\end{document}